\documentclass[12pt]{article}
\usepackage{tikz}
\usetikzlibrary{arrows}
\usepackage[framemethod=tikz]{mdframed}
\usepackage{wrapfig}
\usepackage{amsmath,amssymb}
\usepackage{type1cm}
\usepackage{amsthm}
\usepackage{mathrsfs}
\usepackage{mathtools}
\usepackage{enumerate}
\usepackage[all]{xy} 
\usepackage[vcentermath,enableskew]{youngtab}
\usepackage{ytableau}
\usepackage{diagbox}
\AtBeginDocument{%
   \def\MR#1{}
}
\let\originalthebibliography\thebibliography
\renewcommand{\thebibliography}[1]{%
  \originalthebibliography{#1}%
  \small
  \setlength{\itemsep}{2pt}%
  \setlength{\parsep}{0pt}%
  \setlength{\parskip}{0pt}%
}

\DeclareMathOperator{\ad}{ad}
\DeclareMathOperator{\Ad}{Ad}

\DeclareMathOperator{\en}{end}

\DeclareMathOperator{\supp}{supp}

\newcommand\F{\mathbb{F}}

\newcommand\cB{\mathcal B}
\newcommand\cO{\mathcal O}

\newcommand\cG{\mathcal G}

\newcommand\cT{\mathcal T}

\newcommand\Q{\mathbb Q}

\newcommand\J{\mathbb J}

\newcommand\Z{\mathbb Z}
\newcommand\tW{\tilde W}

\newcommand\tS{\tilde{\mathbb S}}

\newcommand\tp{\mathrm{top}}

\newcommand\unp{\underline p}

\newcommand\vp{\varpi}

\theoremstyle{definition}
\newtheorem{theo}{Theorem}[section]
\newtheorem{prop}[theo]{Proposition}

\newtheorem{lemm}[theo]{Lemma}

\newtheorem{rema}[theo]{Remark}

\newtheorem{thm}{Theorem}[section]

\DeclareMathOperator{\Gal}{Gal}

\newcommand{\scs}{\mathrm{sc}}
\newcommand{\der}{\mathrm{der}}
\newcommand{\reg}{\mathrm{reg}}

\begin{document}
\title{Regular loci of affine Springer fibers at very special parahoric level}
\author{Ryosuke Shimada}
\date{}
\maketitle
\begin{abstract}
For split groups, we prove that centralizers act transitively on regular loci of affine Springer fibers at very special parahoric level.
As an application, we prove Chi's conjecture on the number of top-dimensional irreducible components of affine Lusztig varieties.
\end{abstract}

\section{Introduction}
\label{introduction}

Let $k$ be an algebraically closed field.
We take $L=k((\vp))$, or $L=\operatorname{Frac}(\mathbf{W}(k))$ if $\operatorname{char}k>0$, where $\mathbf{W}(k)$ is the ring of Witt vectors.
We write $\cO$ for the valuation ring of $L$.
Let $\vp$ be a uniformizer of $L$, and let $\cG$ be a split connected reductive group scheme over $\cO$.
We fix a split maximal torus $\cT\subset\cG$ and a Borel subgroup $\cB\subset\cG$ containing $\cT$, both defined over $\cO$.
Let $G,T,B$ denote the generic fibers of $\cG,\cT,\cB$, respectively.
If $\operatorname{char}k>0$, we assume throughout that the order of the Weyl group of $G$ is coprime to $\operatorname{char}k$.
Put $K=\cG(\cO)$.
We choose the Iwahori subgroup $I\subset K$ as in \cite[\S1.1]{He25}.

Let $N$ be the normalizer of $T$ in $G$. The Iwahori--Weyl group is $\tW=N(L)/\cT(\cO)$.
The Iwahori subgroup $I$ determines the length function on $\tW$ and the set $\tS$ of affine simple reflections in $\tW$.
Let $X_*(T)$ denote the cocharacter lattice of $T$.
For $\lambda\in X_*(T)$, we write $\vp^\lambda$ for $\lambda(\vp)\in T(L)$ and for its image in $\tW$.
Let $\tau$ be a length-zero element of $\tW$. We also write $\tau$ for a chosen representative in $N(L)$.

Let $J\subset\tS$ satisfy $\tau J\tau^{-1}=J$ and suppose that $W_J=\langle J\rangle$ is finite.
Let $P\subset G(L)$ be the standard parahoric subgroup generated by $I$ and $W_J$.
Following \cite[\S3.1]{He25}, we call $P$ (or its type $J$) {\it very special} if the longest element of $W_J$ has maximal length among the longest elements of the $W_{J'}$, where $J'$ ranges over all $\tau$-stable subsets of $\tS$ such that $W_{J'}$ is finite.
Throughout the paper, we assume that $P$ is very special.
Let $P^+$ be the pro-unipotent radical of $P$, and let $\overline P=P/P^+$ be its reductive quotient.
Let $\bar\tau$ denote the automorphism of $\overline P$ induced by conjugation by $\tau$.
For $p\in P$, write $\bar p=pP^+$ for its reduction.

For a regular semisimple element $\gamma\in G(L)$, we define the {\it affine Springer fiber} $Y_\tau(\gamma)$ in the partial affine flag variety $G(L)/P$ by
$$
Y_\tau(\gamma)=\{xP\in G(L)/P\mid x^{-1}\gamma x\in P\tau\}.
$$
Its regular locus is defined as in \cite[\S2.4]{He25} by
$$
Y_\tau^{\reg}(\gamma)
=\{xP\in Y_\tau(\gamma)\mid \overline{x^{-1}\gamma x\tau^{-1}}\text{ is $\bar\tau$-regular in }\overline P\}.
$$
We refer to \cite[\S5.5]{HeALV} for the definition of $\bar\tau$-regular elements (cf.\ \S\ref{twisted-regular-subsection}).
The centralizer $Z_G(\gamma)(L)$ acts by left multiplication.

For $P=K$ and $\tau=1$, the following is \cite[Corollary 3.8.2]{Chi22} and \cite[Corollary 7.6]{Chi26}.
See also \cite[Lemma 3.3.1]{Ngo10} and \cite[Theorem 3.5]{He25}.

\begin{thm}[Theorem \ref{parahoric-conjugacy}]\label{main-transitivity}
Let $\gamma\in P\tau$ be regular semisimple.
Then $Y_\tau^{\reg}(\gamma)$ is a single $Z_G(\gamma)(L)$-orbit.
\end{thm}

We reduce the proof of Theorem \ref{main-transitivity} to the analogous statement over $k$ (cf.\ Proposition \ref{finite-lemma}).
For this, we choose a finite tame extension $E/L$ over which the barycenter of the facet of $P$ becomes a hyperspecial vertex.

As an application, we prove Chi's conjecture.
For a dominant cocharacter $\mu\in X_*(T)$ and a regular semisimple element $\gamma\in G(L)$, the {\it affine Lusztig variety} in the affine Grassmannian is
$$
Y_\mu(\gamma)=
\{xK\in G(L)/K\mid x^{-1}\gamma x\in K\vp^\mu K\}.
$$
We regard the locally closed subset $Y_\mu(\gamma)$ as a reduced subscheme in equal characteristic and as a perfect subscheme in mixed characteristic.
The centralizer $Z_G(\gamma)(L)$ acts on $Y_\mu(\gamma)$ by left multiplication.
The variety $Y_\mu(\gamma)$ was first studied by Kottwitz and Viehmann \cite{KV12}.
It naturally arises in the study of orbital integrals of spherical Hecke algebras and plays a role in the conjectural theory of affine character sheaves.

Let $\widehat G$ be the complex Langlands dual group of $G$, and let $V_\mu$ be its irreducible representation of highest weight $\mu$.
For $\lambda\in X_*(T)$, we write $V_\mu(\lambda)$ for the $\lambda$-weight space of $V_\mu$.
For a variety $Y$, we write $\Sigma^{\tp}(Y)$ for the set of its irreducible components of dimension $\dim Y$.

\begin{thm}[Theorem \ref{component-count}]\label{main}
Let $\gamma\in G(L)$ be regular semisimple, and let $\mu\in X_*(T)$ be a dominant cocharacter such that $Y_\mu(\gamma)\ne\varnothing$.
Then
$$
\sharp\bigl(Z_G(\gamma)^\circ(L)\backslash\Sigma^{\tp}(Y_\mu(\gamma))\bigr)
=\dim V_\mu(\lambda_\gamma),
$$
where $\lambda_\gamma$ is a certain weight determined by $\gamma$.
\end{thm}

This theorem was conjectured by Chi \cite[Conjecture 3.9.8]{Chi22}. Chi \cite[Corollary 3.5.3]{Chi22} proved it in equal characteristic when $\gamma$ is conjugate to an element of $T(L)$, and He \cite[Theorem 3.7]{He25} proved it when $\gamma$ has an integral Newton point.

The analogue of Theorem \ref{main} for {\it affine Deligne--Lusztig varieties}, known as the Chen--Zhu conjecture, was proved in \cite{Nie22} and \cite{ZZ20} (cf.\ \S\ref{adlv-subsection}).
He compared affine Lusztig varieties with affine Deligne--Lusztig varieties in \cite{HeALV} and \cite{He25}.
Following the proof of \cite[Theorem 3.7]{He25}, we combine this comparison with Theorem \ref{main-transitivity} to prove Theorem \ref{main}.
When $\gamma$ has an integral Newton point, only the case $P\tau=K$ is needed.

\textbf{Acknowledgments:}
The author would like to thank Kaitao Xie for helpful discussions.

This work was supported by JSPS KAKENHI Grant number JP25K23334.
This paper was written during a stay at UC Berkeley, which was supported by the JSPS Overseas Research Fellowship.
The author would like to thank the university and the host Sug Woo Shin for their hospitality.

OpenAI's ChatGPT was used to assist with mathematical discussion, language, and bibliographic searches.

\section{Proof of Theorem \ref{main-transitivity}}
\label{integral-conjugacy}

\subsection{Twisted regular elements}
\label{twisted-regular-subsection}

Let $k$ be an algebraically closed field, and let $H$ be a simply connected semisimple group over $k$.
For $h\in H$, we write $Z(h)$ for its centralizer in $H$.
Let $T_H$ be a maximal torus of $H$, and write $N_H=\{h\in H\mid hT_H h^{-1}=T_H\}$ for its normalizer.
Let $W_0:=N_H/T_H$ be the Weyl group of $H$ with respect to $T_H$.
If $\operatorname{char}k>0$, assume that $\sharp W_0$ is coprime to $\operatorname{char}k$.
For any $w\in W_0$, we fix a representative $\dot w\in N_H$.
Let $t\in T_H$ be an element of finite order.

Fix $w\in W_0$ such that $t\dot w t^{-1}\dot w^{-1}$ is central in $H$ (we will use a particular choice of $w$ in \S\ref{parahoric-subsection}).
Conjugation by $\dot w$ induces an automorphism $\delta$ of $Z(t)$ given by $\delta(h)=\dot w h\dot w^{-1}$ for $h\in Z(t)$.
Since $H$ is a simply connected semisimple group and $t\in T_H$, the centralizer $Z(t)$ is connected and reductive.
For $h,z\in Z(t)$, the $\delta$-conjugation action is given by $h\cdot_{\delta}z=hz\delta(h^{-1})$.
The $\delta$-centralizer of $z$ is its stabilizer under this action, namely
$$
\{h\in Z(t)\mid h\cdot_{\delta}z=z\}=Z(t)\cap Z(z\dot w).
$$
An element $z\in Z(t)$ is called {\it $\delta$-regular} if its $\delta$-centralizer is of minimal dimension among all the elements of $Z(t)$.
A $\delta$-conjugacy class is regular if it contains a $\delta$-regular element (cf.\ \cite[\S1]{HL26}).

Fix a Borel subgroup $B_H$ of $H$ containing $T_H$. An element of $W_0$ is called a {\it Coxeter element} if it is a product of the simple reflections determined by $B_H$, each occurring exactly once.
Put $B_t=Z(t)\cap B_H$.
Then $B_t$ is a Borel subgroup of $Z(t)$ containing $T_H$.
We further assume that $\delta(B_t)=\dot w B_t\dot w^{-1}=B_t$.

\begin{prop}\label{finite-lemma}
Suppose that $\dot c\dot w^{-1}\in Z(t)$ for a Coxeter element $c\in W_0$.
Then the following statements hold.
\begin{enumerate}[(i)]
\item If $z\in Z(t)$ is $\delta$-regular, then $z\dot w$ is regular in $H$.
\item If $z_1,z_2\in Z(t)$ are $\delta$-regular and $z_1\dot w,z_2\dot w$ are conjugate in $H$, then there exists $h\in Z(t)$ such that $h(z_1\dot w)h^{-1}=z_2\dot w$.
\end{enumerate}
\end{prop}

\begin{proof}
Let $z\in Z(t)$ be $\delta$-regular. Then by \cite[\S5.5(a)]{HeALV}, the Lusztig variety
$$
\{xB_t\in Z(t)/B_t\mid x^{-1}z\delta(x)\in B_t\dot c\dot w^{-1}B_t\}
$$
is nonempty.
Hence there exists $h\in Z(t)$ such that $h\cdot_\delta z\in B_t\dot c\dot w^{-1}B_t$.
Since $\delta(B_t)=B_t$, we have $h(z\dot w)h^{-1}=(h\cdot_\delta z)\dot w\in B_t\dot c B_t\subset B_H\dot c B_H$.
Let $U$ be the unipotent radical of $B_H$. Since the map $T_H\to T_H$, $s\mapsto s\dot c s^{-1}\dot c^{-1}$, is surjective by \cite[Lemma 7.6]{Steinberg65}, we may further conjugate by $T_H\subset Z(t)$ to arrange that $z\dot w\in U\dot cU$.

Let $U^-$ be the unipotent radical of the opposite Borel subgroup. Put $\operatorname{St}=(U\cap\dot cU^-\dot c^{-1})\dot c$.
By \cite[Proposition 8.9]{Steinberg65}, there is a unique pair $(u,a)\in U\times\operatorname{St}$ such that $z\dot w=uau^{-1}$.
By \cite[Theorem 1.4]{Steinberg65}, $a$ is regular in $H$, which proves (i).

Put $\zeta=t\dot w t^{-1}\dot w^{-1}$, which is central in $H$. Since $z\in Z(t)$, we have $t(z\dot w)t^{-1}=\zeta z\dot w$. We conjugate $z\dot w=uau^{-1}$ by $t$ and multiply by $\zeta^{-1}$ on the left to obtain
$$
z\dot w=(tut^{-1})(\zeta^{-1}tat^{-1})(tut^{-1})^{-1}.
$$
Since $t\in T_H$, we have $tut^{-1}\in U$. The assumption $\dot c\dot w^{-1}\in Z(t)$ gives $t\dot c t^{-1}=\zeta\dot c$.
Thus $\zeta^{-1}tat^{-1}=t(a\dot c^{-1})t^{-1}\dot c\in\operatorname{St}$, since $T_H$ normalizes $U\cap\dot cU^-\dot c^{-1}$. By uniqueness of the pair $(u,a)$, we have $u=tut^{-1}$ and $a=\zeta^{-1}tat^{-1}$. Hence $u\in Z(t)$. Thus $z\dot w$ is conjugate under $Z(t)$ to $a\in\operatorname{St}$.

Let $z_1,z_2\in Z(t)$ be two $\delta$-regular elements. We have $z_i\dot w=u_ia_iu_i^{-1}$ for some $u_i\in Z(t)$ and $a_i\in\operatorname{St}$.
Since $z_1\dot w$ and $z_2\dot w$ are conjugate in $H$, so are $a_1$ and $a_2$.
Again by \cite[Theorem 1.4]{Steinberg65}, their common conjugacy class meets $\operatorname{St}$ in only one point since $H$ is semisimple and simply connected. Hence $a_1=a_2$.
Thus $(u_2u_1^{-1})\cdot_\delta z_1=z_2$.
\end{proof}

\subsection{Regular locus at very special parahoric level}
\label{parahoric-subsection}

We keep the notation of \S\ref{introduction}.
We write $G^{\scs}\longrightarrow G^{\der}$ for the simply connected covering of the derived subgroup of $G$.
Let $T^{\scs}$ be the inverse image of $T\cap G^{\der}$ in $G^{\scs}$.
By abuse of notation, we shall often use the same letter to denote an element in $G^{\scs}$ and its image in $G$.
On the other hand, we use the superscript $\ad$ to denote images under the natural maps $G^{\scs}\to G^{\ad}$ and $G\to G^{\ad}$.

Let $\cG^{\scs}$ and $\cG^{\ad}$ be the split hyperspecial models of $G^{\scs}$ and $G^{\ad}$ over $\cO$ corresponding to the vertex of $K$ in the Bruhat--Tits building of $G$ over $L$.
Let $\mathcal T^{\scs}\subset\cG^{\scs}$ be the split maximal torus with generic fiber $T^{\scs}$.
Let $H$ and $T_H$ be the special fibers of $\cG^{\scs}$ and $\mathcal T^{\scs}$, respectively.
Then $H$ is semisimple and simply connected with maximal torus $T_H$.
For any finite extension $E/L$ with valuation ring $\cO_E$ and any $g\in \cG^{\scs}(\cO_E)$, we write $\bar g\in H(k)$ for its reduction.
We identify $N_H/T_H$ with $N/T$ and denote both Weyl groups by $W_0$.
We write $\tW=X_*(T)\rtimes W_0$.

The goal of this subsection is to prove the following theorem.

\begin{theo}\label{parahoric-conjugacy}
Let $\gamma_1,\gamma_2\in P\tau$ be regular semisimple elements such that $\overline{\gamma_1\tau^{-1}}$ and $\overline{\gamma_2\tau^{-1}}$ are $\bar\tau$-regular in $\overline P$.
If $\gamma_2=g\gamma_1g^{-1}$ for some $g\in G(L)$, then there exists $p\in P$ such that
$$
\gamma_2=p\gamma_1p^{-1}.
$$
As a consequence, if $\gamma\in P\tau$ is regular semisimple, then the set of closed points of $Y_\tau^{\reg}(\gamma)$ is a single $Z_G(\gamma)(L)$-orbit.
\end{theo}

To find $p\in P$ with $p\gamma_1p^{-1}=\gamma_2$, we first choose a tamely ramified extension $E/L$ and $r\in T^{\scs}(E)$ such that $(r\gamma_1r^{-1})^{\ad},(r\gamma_2r^{-1})^{\ad}\in \cG^{\ad}(\cO_E)$.
We lift these images to $g_1,g_2\in \cG^{\scs}(\cO_E)$.
We use Proposition \ref{finite-lemma} to find $h\in H(k)$ with $h\bar g_1h^{-1}=\bar g_2$.
We then lift $h$ to $\widetilde h\in \cG^{\scs}(\cO_E)$ with $\widetilde h g_1\widetilde h^{-1}=g_2$ and $r^{-1}\widetilde h r\in G^{\scs}(L)$.
The image of $r^{-1}\widetilde h r$ in $G(L)$ is the required conjugating element in $P$.
We give the details below.

Let $v$ be the barycenter of the facet of $P$ in the apartment corresponding to the image $T^{\ad}$ of $T$ in $G^{\ad}$.
Since $P$ is $\tau$-stable, we have $\tau v=v$.
We designate the vertex of $K$ as the origin. This identifies the apartment with $X_*(T^{\ad})_{\mathbb R}$.
We use the simple coroots determined by $B$ as a basis.
\begin{lemm}\label{tame-extension}
There exist a finite tamely ramified extension $E/L$ and $r\in T^{\scs}(E)$ such that $r\cdot v=0$.
\end{lemm}

A similar argument can be found in the proof of \cite[Proposition 2.8]{PR24}.

\begin{proof}
We may assume that $G$ is adjoint.
Assume further that $G$ is simple.
The vertices $v_1,\ldots,v_m$ of the facet of $P$ form one $\tau$-orbit since $P$ is very special, so $m$ divides $\sharp W_0$.
The stabilizer $W_{v_i}$ of $v_i$ in the affine Weyl group maps injectively to $W_0$ under the natural projection.
Thus $\sharp W_{v_i}$ divides $\sharp W_0$, and
$$
(\sharp W_{v_i})v_i=\sum_{u\in W_{v_i}}u(0)\in X_*(T^{\scs}).
$$
Hence $ev\in X_*(T^{\scs})$ for $e=(\sharp W_0)^2$, since $v=m^{-1}\sum_i v_i$.
Choose $\vp_E^e=\vp$ and set $E=L(\vp_E)$ and $r=\vp_E^{-ev}$. Then we have $r\cdot v=0$.
Since $\operatorname{char}k\nmid \sharp W_0$, $E/L$ is a tamely ramified extension.

If $G$ is not simple, we decompose $G$ into its simple factors. We then decompose $G^{\scs}$, $T$, $T^{\scs}$ and the apartment accordingly.
The image of $P$ in each factor of $G(L)$ is a very special parahoric subgroup.
The composite of tamely ramified extensions is tamely ramified.
For each simple factor, we choose a tame extension and a torus element as above.
We take $E$ to be the composite of these extensions and $r$ to be the product of these torus elements.
This completes the proof.
\end{proof}

Choose $E/L$ and $r$ as in Lemma \ref{tame-extension}. For any $\gamma\in P\tau$, the adjoint image $(r\gamma r^{-1})^{\ad}$ fixes the origin and belongs to $\cG^{\ad}(\cO_E)$.
In particular, we may choose a representative $\dot w\in N_H$ of the finite part $w\in W_0$ of $\tau$ lifting $\overline{(r\tau r^{-1})^{\ad}}$.
Note that $H^{\ad}=\cG^{\ad}_k$.

Since $k$ is algebraically closed, the tamely ramified extension $E/L$ is a cyclic Galois extension. Choose a generator $\eta$ of $\Gal(E/L)$ and put $t=r\eta(r)^{-1}\in \mathcal T^{\scs}(\cO_E)$.
By definition, $\bar t\dot w \bar t^{-1}\dot w^{-1}$ is central in $H$.
The element $\bar t\in T_H$ has finite order. Indeed,
$$
\bar t^{[E:L]}=\prod_{j=0}^{[E:L]-1}\overline{\eta^j(t)}=\prod_{j=0}^{[E:L]-1}\overline{\eta^j(r)\eta^{j+1}(r)^{-1}}=\overline{r\eta^{[E:L]}(r)^{-1}}=\overline{rr^{-1}}=1.
$$

Define the automorphism $F\colon G^{\scs}(E)\longrightarrow G^{\scs}(E)$ by
$$
F(g)=r\eta(r^{-1}gr)r^{-1}=t\eta(g)t^{-1}.
$$
We define $F$ on $G(E)$ and $G^{\ad}(E)$ in the same way.
Let $\mathcal P^{\scs}$ be the parahoric model of $G^{\scs}$ over $\cO$ at $v$, so that $\mathcal P^{\scs}(\cO)=G^{\scs}(L)\cap r^{-1}\cG^{\scs}(\cO_E)r$.
Put $\mathcal R=\operatorname{Res}_{\cO_E/\cO}(\cG^{\scs}_{\cO_E})$.
Let $\Gamma=\Gal(E/L)$ act on $\mathcal R$, with $\eta$ acting by $F$.
\begin{prop}\label{parahoric-descent}
The isomorphism $\mathcal P^{\scs}_L=G^{\scs}\xrightarrow{\sim}\mathcal R_L^\Gamma$ given by $g\mapsto rgr^{-1}$ extends to an isomorphism $\mathcal P^{\scs}\xrightarrow{\sim}\mathcal R^\Gamma$ of $\cO$-group schemes.
As a consequence, we have
$$
\mathcal P^{\scs}_k/R_u(\mathcal P^{\scs}_k)\cong Z(\bar t).
$$
\end{prop}

The first assertion can be found in \cite[Proposition 2.8]{PR24}.

\begin{proof}
Since $\sharp\Gamma$ is invertible in $k$, $\mathcal R^\Gamma$ is smooth over $\cO$ by \cite[Proposition 3.4]{Edi92}.
Thus the first assertion follows from $r\mathcal P^{\scs}(\cO)r^{-1}=\cG^{\scs}(\cO_E)^F=\mathcal R^\Gamma(\cO)$.
By \cite[Proposition 9.11.1(3)]{KP23}, we have $\mathcal R_k/R_u(\mathcal R_k)\cong H$.
Under this isomorphism, $\eta$ acts on $H$ by $\Ad(\bar t)$.
Since $\sharp\Gamma$ is invertible in $k$, we obtain
$$
\mathcal P^{\scs}_k/R_u(\mathcal P^{\scs}_k)\cong H^\Gamma=Z(\bar t)
$$
by \cite[Lemma 12.5.4(1)]{KP23}.
\end{proof}

By Proposition \ref{parahoric-descent}, we obtain the following identification of root systems
$$
\Phi(Z(\bar t),T_H)=\{\alpha\in\Phi(H,T_H)\mid\langle\alpha,v\rangle\in\Z\}=\{\alpha\in\Phi(G,T)\mid\langle\alpha,v\rangle\in\Z\}.
$$
Moreover, the central isogeny $G^{\scs}\times Z(G)^\circ\to G$ induces a central isogeny
$$
Z(\bar t)\times Z(\cG)^\circ_k\longrightarrow\overline P
$$
compatible with $\delta\times\mathrm{id}$ and $\bar\tau$, where $\delta(h)=\dot w h\dot w^{-1}$ for $h\in Z(\bar t)$.
Here $Z(G)$ and $Z(\cG)$ denote the centers of $G$ and $\cG$, respectively.
Let $B_{\bar t}\supset T_H$ denote the Borel subgroup of $Z(\bar t)$ determined by the inverse image of $I/P^+$.
Since $\tau I\tau^{-1}=I$, we have $\delta(B_{\bar t})=\dot w B_{\bar t}\dot w^{-1}=B_{\bar t}$.

Write $\tW=W_a\rtimes\Omega$, with affine Weyl group $W_a=\langle\tS\rangle$ and subgroup $\Omega$ of length-zero elements.
The proof of Lemma \ref{parahoric-coxeter} relies on the following result, which follows immediately from \cite[Proposition 3.4.4]{HZZ25}.

\begin{prop}\label{finite-coxeter-part}
Suppose that $G$ is adjoint and simple. Let $u\tau\in W_a\tau$ be an element whose finite part is a Coxeter element of $W_0$. Then the $W_a$-conjugacy class of $u\tau$ meets $W_{J_0}\tau$ for some very special type $J_0\subset\tS$.
\end{prop}

\begin{lemm}\label{parahoric-coxeter}
There exist a Borel subgroup $B_H$ of $H$ containing $T_H$ and a Coxeter element $c\in W_0$ with respect to $B_H$ such that $Z(\bar t)\cap B_H=B_{\bar t}$ and $\dot c\dot w^{-1}\in Z(\bar t)$.
\end{lemm}

\begin{proof}
We may assume that $G$ is adjoint and simple. Choose a Coxeter element $c_0$ with respect to $B$. The element $c_0w^{-1}\tau\in W_a\tau$ has finite part $c_0$. By Proposition \ref{finite-coxeter-part}, its $W_a$-conjugacy class meets $W_{J_0}\tau$ for some very special type $J_0$.
As noted at the end of \cite[\S1.1]{HZZ25}, we have $J=\omega J_0\omega^{-1}$ for some length-zero element $\omega\in\Omega$. Hence $W_J\tau=\omega(W_{J_0}\tau)\omega^{-1}$ contains an element $u\tau$ whose finite part $c'$ is conjugate to $c_0$.
Then $c'w^{-1}$, the finite part of $u$, lies in the Weyl group of $Z(\bar t)$.

Choose a Borel subgroup $B'_H$ containing $T_H$ for which $c'$ is a Coxeter element. Since $Z(\bar t)\cap B'_H$ and $B_{\bar t}$ are Borels containing $T_H$, there is an $n\in Z(\bar t)\cap N_H$ with $n(Z(\bar t)\cap B'_H)n^{-1}=B_{\bar t}$. Set $B_H=nB'_Hn^{-1}$ and let $c$ be the image of $n\dot c'n^{-1}$ in $W_0$. Since $w$ normalizes $W(Z(\bar t),T_H)$, we have $cw^{-1}\in W(Z(\bar t),T_H)$, or equivalently $\dot c\dot w^{-1}\in Z(\bar t)$. By construction, $c$ is a Coxeter element with respect to $B_H$.
\end{proof}

We choose $B_H$ and $c$ as in Lemma \ref{parahoric-coxeter}, so that $\bar t,w,B_H,c$ satisfy the hypotheses of Proposition \ref{finite-lemma}.

\begin{proof}[Proof of Theorem \ref{parahoric-conjugacy}]
Choose $g\in G(L)$ with $\gamma_2=g\gamma_1g^{-1}$. By our assumption on $\operatorname{char}k$, the covering $\cG^{\scs}_{\cO_E}\to \cG^{\ad}_{\cO_E}$ is finite \'etale.
We may therefore lift $(r\gamma_1r^{-1})^{\ad}$ to $g_1\in \cG^{\scs}(\cO_E)$. Define $g_2=\operatorname{Ad}((rgr^{-1})^{\ad})(g_1)\in G^{\scs}(E)$. Then $g_2^{\ad}=(r\gamma_2r^{-1})^{\ad}\in \cG^{\ad}(\cO_E)$, so $g_2\in \cG^{\scs}(\cO_E)$ by properness of the covering.

Put $z_i=\bar g_i\dot w^{-1}$. We have $z_i^{\ad}=\overline{(r\gamma_i\tau^{-1}r^{-1})^{\ad}}\in Z(\bar t)/Z(H)$ by Proposition \ref{parahoric-descent}. Thus $z_i\in Z(\bar t)$. Since $\overline{\gamma_i\tau^{-1}}$ is $\bar\tau$-regular, $z_i$ is $\delta$-regular.
By Proposition \ref{finite-lemma}(i), $\bar g_i=z_i\dot w$ is regular in $H$. Let $\chi:\cG^{\scs}\to\cG^{\scs}/\!/\cG^{\scs}$ be the Steinberg morphism over $\cO$. Since $g_1$ and $g_2$ are geometrically conjugate, we have $\chi(g_1)=\chi(g_2)$. By \cite[Theorem 4.1(1)]{Lee15}, $(\cG^{\scs}/\!/\cG^{\scs})_k$ is canonically isomorphic to $H/\!/H$, and under this identification the base change $\chi_k$ of $\chi$ to $k$ is the Steinberg morphism for $H$. Thus $\chi_k(\bar g_1)=\overline{\chi(g_1)}=\overline{\chi(g_2)}=\chi_k(\bar g_2)$. Hence $\bar g_1,\bar g_2$ are conjugate in $H$ by \cite[Theorem 6.11(c)]{Steinberg65}.
By Proposition \ref{finite-lemma}(ii), we may choose $h\in Z(\bar t)$ such that $h\bar g_1h^{-1}=\bar g_2$.

Let $(\cG^{\scs}_{\cO_E})^{\reg}$ denote the regular locus of $\cG^{\scs}_{\cO_E}$.
Since $\bar g_i$ is regular, the section $g_i\colon\operatorname{Spec}\cO_E\to\cG^{\scs}_{\cO_E}$ factors through $(\cG^{\scs}_{\cO_E})^{\reg}$.
By the proof of \cite[Theorem 3.7]{Shotton24}, the conjugation morphism
$$
c\colon \cG^{\scs}_{\cO_E}\times_{\cO_E}(\cG^{\scs}_{\cO_E})^{\reg}
\longrightarrow
(\cG^{\scs}_{\cO_E})^{\reg}\times_{\cG^{\scs}_{\cO_E}/\!/\cG^{\scs}_{\cO_E}}(\cG^{\scs}_{\cO_E})^{\reg},
\qquad (x,g)\longmapsto (g,xgx^{-1}),
$$
is smooth.
Since $\chi(g_1)=\chi(g_2)$, the pair $(g_1,g_2)$ gives a section of the target of $c$.
As $c$ is smooth, so is its pullback $\mathcal X\to\operatorname{Spec}\cO_E$ along this section.
The pullback of $g_2$ along the orbit map $\cG^{\scs}_{\cO_E}\to\cG^{\scs}_{\cO_E}$, $x\mapsto xg_1x^{-1}$, identifies $\mathcal X$ with a closed subscheme of $\cG^{\scs}_{\cO_E}$.
We have
$$
\mathcal X(\cO_E)=\{x\in\cG^{\scs}(\cO_E)\mid xg_1x^{-1}=g_2\},\quad \mathcal X(k)=\{x\in H(k)\mid x\bar g_1x^{-1}=\bar g_2\}.
$$

As $F$ fixes $r\gamma_1r^{-1}$ and $rgr^{-1}$, we have $F(g_1)g_1^{-1}=F(g_2)g_2^{-1}\in Z(G^{\scs})(E)$. Thus $\mathcal X$ is stable under $F$. The point $h\in\mathcal X(k)$ satisfies $F(h)=\bar t h\bar t^{-1}=h$ since $h\in Z(\bar t)$. By \cite[Proposition 3.3]{Hartmann12}, $h$ lifts to an element $\widetilde h\in\mathcal X(\cO_E)$ with $F(\widetilde h)=\widetilde h$.
Put $p=r^{-1}\widetilde h r$. Then $\eta(p)=r^{-1}F(\widetilde h)r=p$, so $p\in G^{\scs}(L)\cap r^{-1}\cG^{\scs}(\cO_E)r=\mathcal P^{\scs}(\cO)$. Hence its image in $G(L)$ lies in $P$. By the choice of $g_2$, we have $(r\gamma_1r^{-1})g_1^{-1}=(r\gamma_2r^{-1})g_2^{-1}\in Z(G)(E)$. Thus $\widetilde h g_1\widetilde h^{-1}=g_2$ gives $p\gamma_1p^{-1}=\gamma_2$, as desired.

Let $\gamma\in P\tau$ be regular semisimple. The regular locus is nonempty by \cite[Theorem 3.1(2)]{HeALV}.
Let $x_1P,x_2P\in Y_\tau^{\reg}(\gamma)\subset G(L)/P$ and put $\gamma_i=x_i^{-1}\gamma x_i\in P\tau$. By definition, $\overline{\gamma_i\tau^{-1}}$ is $\bar\tau$-regular.
Since $x_2^{-1}x_1$ conjugates $\gamma_1$ to $\gamma_2$, the first assertion gives $p\in P$ with $p\gamma_1p^{-1}=\gamma_2$.
Then $x_2px_1^{-1}\in Z_G(\gamma)(L)$ carries $x_1P$ to $x_2P$.
This completes the proof.
\end{proof}

\begin{rema}
If $P=K$ and $\tau=1$, we may take $E=L$ and $r=t=\dot w=1$, so that $Z(\bar t)=H$ and $\delta=\mathrm{id}$.
In this case, Proposition \ref{finite-lemma} is trivial, and our proof is similar to Chi's proof of \cite[Lemma 7.5]{Chi26}.
\end{rema}

\section{Proof of Theorem \ref{main}}
\label{proof-main}

\subsection{Affine Deligne--Lusztig varieties}
\label{adlv-subsection}

Let $\cG'$ be a split connected reductive group scheme over $\F_q[[\vp]]$ with generic fiber $G'$.
We fix a split maximal torus $\cT'\subset\cG'$ and a Borel subgroup $\cB'\subset\cG'$ containing $\cT'$.
Let $T'$ and $B'$ denote their generic fibers, respectively.
Put $L'=\overline{\F}_q((\vp))$ and let $\sigma$ denote the Frobenius automorphism of $L'$.
Fix $b\in G'(L')$ and write $[b]$ for its $\sigma$-conjugacy class.
Put $K'=\cG'(\overline{\F}_q[[\vp]])$ and choose a $\sigma$-stable Iwahori subgroup $I'\subset K'$.

In this subsection, $\tW$ and $W_0$ denote the Iwahori--Weyl group and the finite Weyl group of $G'$, respectively.
For $w\in\tW$ and dominant $\mu\in X_*(T')$, we define the affine Deligne--Lusztig varieties at Iwahori and hyperspecial level by
$$
\begin{aligned}
X_w(b)&=\{xI'\in G'(L')/I'\mid x^{-1}b\sigma(x)\in I'wI'\},\\
X_\mu(b)&=\{xK'\in G'(L')/K'\mid x^{-1}b\sigma(x)\in K'\vp^\mu K'\}.
\end{aligned}
$$
The group $\J_b=\{j\in G'(L')\mid j^{-1}b\sigma(j)=b\}$ acts on both by left multiplication.

Let $\kappa(b)\in\pi_1(G')$ and $\nu_b\in X_*(T')_{\Q}$ be the Kottwitz invariant and the dominant Newton point of $[b]$.
For $\lambda,\nu\in X_*(T')_{\Q}$, write $\lambda\le\nu$ if $\nu-\lambda$ is a nonnegative rational linear combination of simple coroots.
By Rapoport--Richartz \cite[Theorem 4.2]{RR96} and Gashi \cite[Theorem 1.1]{Gashi10}, $X_\mu(b)\ne\varnothing$ if and only if $\kappa(\vp^\mu)=\kappa(b)$ and $\nu_b\le\mu$.

For $w\in\tW$, we fix a reduction tree $\mathcal T$ of $w$.
The vertices of $\mathcal T$ are elements of $\tW$, and its directed edges are of the form $x\rightharpoonup y$.
For the precise definition of a reduction tree, see \cite[\S2.2]{He25}.
A reduction path in $\mathcal T$ is a path $\unp:w\rightharpoonup w_1\rightharpoonup\cdots\rightharpoonup w_n$ from $w$ to a vertex $w_n$ with no edge of the form $w_n\rightharpoonup y$.
We write $\en(\unp)=w_n$ and $\ell(\unp)=n$.
The element $\en(\unp)$ has minimal length in its conjugacy class.
Let $\mathrm{RP}(\mathcal T)$ denote the set of reduction paths in $\mathcal T$.

For $x\in\tW$, choose $m>0$ such that $x^m=\vp^\lambda$ for some $\lambda\in X_*(T')$, and put $\nu_x=\lambda/m$.
Let $\bar\nu_x$ be the dominant representative of $\nu_x$.
The Kottwitz map on $\tW$ is defined by $\kappa(x):=xW_a\in\tW/W_a\cong\pi_1(G')$, where $W_a$ is the affine Weyl group.
Let $\mathrm{RP}_{[b]}(\mathcal T)$ be the set of paths $\unp\in\mathrm{RP}(\mathcal T)$ such that
$$
\bigl(\kappa(\en(\unp)),\bar\nu_{\en(\unp)}\bigr)=(\kappa(b),\nu_b).
$$

\begin{prop}\label{adlv-grassmannian-path-count}
Let $w_0$ be the longest element of $W_0$.
Let $\mu\in X_*(T')$ be dominant with $X_\mu(b)\ne\varnothing$, and put $w=w_0\vp^\mu$.
Then
$$
\sharp\bigl(\J_b\backslash\Sigma^{\tp}(X_\mu(b))\bigr)
=\sharp\{\unp\in\mathrm{RP}_{[b]}(\mathcal T)\mid
\ell(\unp)+\ell(\en(\unp))-\langle\bar\nu_{\en(\unp)},2\rho\rangle=\dim X_w(b)\},
$$
where $\rho$ is the half-sum of the positive roots.
\end{prop}

\begin{proof}
By \cite[Proposition 2.4.10]{HZZ25}, the projection $xI'\mapsto xK'$ induces a $\J_b$-equivariant bijection
$$
\Sigma^{\tp}(X_w(b))\xrightarrow{\sim}\Sigma^{\tp}(X_\mu(b)).
$$
Thus the formula follows from \cite[Proposition 2.5]{He25}.
\end{proof}

By \cite[Lemma/Definition 2.1]{HV18}, the set
$$
\{\lambda\in X_*(T')\mid \kappa(\vp^\lambda)=\kappa(b),\ \lambda\le\nu_b\}
$$
has a unique maximum $\lambda_b$.
Let $V_\mu^{\widehat{G'}}$ be the irreducible representation of the complex dual group $\widehat{G'}$ with highest weight $\mu$.
By Nie \cite{Nie22} and Zhou--Zhu \cite{ZZ20}, we have
$$
\sharp\bigl(\J_b\backslash\Sigma^{\tp}(X_\mu(b))\bigr)=\dim V_\mu^{\widehat{G'}}(\lambda_b).
$$
This equality was conjectured by Miaofen Chen and Xinwen Zhu.
In particular,
$$
\bigl(\kappa(\vp^\mu)=\kappa(b)\text{ and }\nu_b\le\mu\bigr) \Longleftrightarrow X_\mu(b)\ne\varnothing \Longleftrightarrow V_\mu^{\widehat{G'}}(\lambda_b)\ne 0.
$$

\subsection{Affine Lusztig varieties}

We now take $G'$ with the same based root datum as $G$.
We identify the Iwahori--Weyl groups of $G$ and $G'$ by a length-preserving isomorphism induced by an isomorphism of their based root data.

Fix a regular semisimple element $\gamma\in G(L)$ and write $\{\gamma\}$ for its conjugacy class.
We define the Kottwitz invariant $\kappa(\gamma)\in\pi_1(G)$ and the dominant Newton point $\nu_\gamma\in X_*(T)_{\Q}$ of $\{\gamma\}$ as in \cite{KV12}.
By \cite[Theorem 6.9]{HeALV}, for any dominant cocharacter $\mu\in X_*(T)$, $Y_\mu(\gamma)\ne\varnothing$ if and only if $\kappa(\vp^\mu)=\kappa(\gamma)$ and $\nu_\gamma\le\mu$.
Let $B(G)$ be the set of conjugacy classes in $G(L)$ and $B_\sigma(G')$ the set of $\sigma$-conjugacy classes in $G'(L')$.
The map $f_{G,G'}\colon B(G)\longrightarrow B_\sigma(G')$ of \cite[\S1.2]{He25} sends $\{\gamma\}$ to the unique $\sigma$-conjugacy class $[b]$ with $\kappa(b)=\kappa(\gamma)$ and $\nu_b=\nu_\gamma$.
Choose $b\in G'(L')$ representing this class.

For $w\in\tW$ such that $X_w(b)\ne\varnothing$, let $\mathcal T$ be the reduction tree of $w$ fixed in \S\ref{adlv-subsection}.
For each $\unp\in\mathrm{RP}_{[b]}(\mathcal T)$, let $M_{\unp}$ be the Levi subgroup of $G$ generated by $T$ and the root subgroups $U_\alpha$ with $\langle\alpha,\nu_{\en(\unp)}\rangle=0$.
By \cite[\S4.2]{HeALV}, we have $\en(\unp)=u\tau$, where $u$ belongs to the affine Weyl group of $M_{\unp}$ and $\tau$ is a length-zero element of the Iwahori--Weyl group of $M_{\unp}$.
Let $J_{\unp}=\bigcup_{n\in\mathbb Z}\tau^n\supp(u)\tau^{-n}$, where $\supp(u)$ is the set of simple reflections that occur in some (or equivalently, any) reduced expression for $u$.
Let $P_{\unp}\subset M_{\unp}(L)$ be the standard parahoric subgroup of type $J_{\unp}$.

By \cite[Theorem 4.4]{HeALV}, we may choose a $G(L)$-conjugate $\gamma'$ of $\gamma$ in the $(I\cap M_{\unp}(L))$-double coset of $\en(\unp)$, which is contained in $P_{\unp}\tau$.
We define the affine Springer fiber $Y^{M_{\unp}}_{J_{\unp},\tau}(\gamma')$ and the integer $n_{\unp,\gamma}$ by
$$
\begin{aligned}
Y^{M_{\unp}}_{J_{\unp},\tau}(\gamma')
 &=\{xP_{\unp}\in M_{\unp}(L)/P_{\unp}\mid x^{-1}\gamma'x\in P_{\unp}\tau\},\\
n_{\unp,\gamma}
 &=\sharp\bigl(Z_{M_{\unp}}(\gamma')(L)\backslash
       \Sigma^{\tp}(Y^{M_{\unp}}_{J_{\unp},\tau}(\gamma'))\bigr).
\end{aligned}
$$

For $w\in\tW$, the affine Lusztig variety $Y_w(\gamma)$ in the affine flag variety is defined by
$$
Y_w(\gamma)=\{xI\in G(L)/I\mid x^{-1}\gamma x\in IwI\}.
$$

\begin{prop}\label{alv-grassmannian-path-count}
Let $\mu\in X_*(T)$ be dominant with $Y_\mu(\gamma)\ne\varnothing$, and put $w=w_0\vp^\mu$.
Then $X_\mu(b)\ne\varnothing$ and
$$
\sharp\bigl(Z_G(\gamma)(L)\backslash\Sigma^{\tp}(Y_\mu(\gamma))\bigr)
=\sum_{\substack{\unp\in\mathrm{RP}_{[b]}(\mathcal T)\\
\ell(\unp)+\ell(\en(\unp))-\langle\bar\nu_{\en(\unp)},2\rho\rangle=\dim X_w(b)}}
n_{\unp,\gamma}.
$$
\end{prop}

\begin{proof}
By the proof of \cite[Theorem 3.7]{He25}, the map $xI\mapsto xK$ induces a $Z_G(\gamma)(L)$-equivariant bijection
$$
\Sigma^{\tp}(Y_w(\gamma))\xrightarrow{\sim}\Sigma^{\tp}(Y_\mu(\gamma)).
$$
Thus the result follows from \cite[Theorem 2.8]{He25} and \cite[Theorem 6.7]{HeALV}.
\end{proof}

By the argument of \cite[Lemma 3.9.7]{Chi22}, the set
$$
\{\lambda\in X_*(T)\mid \lambda\text{ is dominant},\ \kappa(\vp^\lambda)=\kappa(\gamma),\ \nu_\gamma\le\lambda\}
$$
has a unique minimum $\lambda_\gamma$.
Let $\bar\lambda_b$ be the dominant representative of $W_0\lambda_b$.
Then $\bar\lambda_b$ is the minimal dominant cocharacter $\lambda$ such that $V_\lambda^{\widehat{G'}}(\lambda_b)\ne 0$.
The equivalences at the end of \S\ref{adlv-subsection} therefore give $\lambda_\gamma=\bar\lambda_b$.
The following is \cite[Conjecture 3.9.8]{Chi22}.
\begin{theo}\label{component-count}
Let $\gamma\in G(L)$ be regular semisimple, and let $\mu\in X_*(T)$ be a dominant cocharacter such that $Y_\mu(\gamma)\ne\varnothing$.
Then
$$
\sharp\bigl(Z_G(\gamma)^\circ(L)\backslash\Sigma^{\tp}(Y_\mu(\gamma))\bigr)
=\sharp\bigl(Z_G(\gamma)(L)\backslash\Sigma^{\tp}(Y_\mu(\gamma))\bigr)
=\dim V_\mu^{\widehat G}(\lambda_\gamma).
$$
\end{theo}

\begin{proof}
We first prove the second equality.
Let $\mathcal T$ be a reduction tree of $w_0\vp^\mu$ and let $\unp\in\mathrm{RP}_{[b]}(\mathcal T)$ satisfy $\ell(\unp)+\ell(\en(\unp))-\langle\bar\nu_{\en(\unp)},2\rho\rangle=\dim X_{w_0\vp^\mu}(b)$.
By \cite[Proposition 3.1]{He25}, $P_{\unp}$ is a very special parahoric subgroup of $M_{\unp}(L)$.
Since the set of regular elements in any reductive group over $k$ is open, the regular locus $Y_{\unp,\tau}^{\reg}(\gamma')$ of $Y^{M_{\unp}}_{J_{\unp},\tau}(\gamma')$ is open.
By the proof of \cite[Proposition 2.7]{He25}, the inclusion induces a natural bijection $\Sigma^{\tp}(Y_{\unp,\tau}^{\reg}(\gamma'))\xrightarrow{\sim}\Sigma^{\tp}(Y^{M_{\unp}}_{J_{\unp},\tau}(\gamma'))$.
Since $\Sigma^{\tp}(Y^{M_{\unp}}_{J_{\unp},\tau}(\gamma'))$ is locally finite (cf.\ \cite[Theorem 6.10]{Chi26}), each $C\in\Sigma^{\tp}(Y^{M_{\unp}}_{J_{\unp},\tau}(\gamma'))$ contains a regular closed point $x$ such that $x\notin C'$ for every $C'\in\Sigma^{\tp}(Y^{M_{\unp}}_{J_{\unp},\tau}(\gamma'))\setminus\{C\}$.
By Theorem \ref{parahoric-conjugacy}, $Y_{\unp,\tau}^{\reg}(\gamma')$ is a single $Z_{M_{\unp}}(\gamma')(L)$-orbit.
Thus $n_{\unp,\gamma}=1$.
The Chen--Zhu conjecture combined with Proposition~\ref{adlv-grassmannian-path-count} and Proposition~\ref{alv-grassmannian-path-count} gives
\begin{equation*}
\sharp\bigl(Z_G(\gamma)(L)\backslash\Sigma^{\tp}(Y_\mu(\gamma))\bigr)
=\dim V_\mu^{\widehat{G'}}(\lambda_b)
=\dim V_\mu^{\widehat G}(\lambda_\gamma).
\end{equation*}

For the first equality, let $\pi\colon\widetilde G\to G$ be a $z$-extension whose kernel is an $L$-split torus $\widetilde T$.
By Hilbert 90, $\pi\colon\widetilde G(L)\to G(L)$ is surjective.
Choose a lift $\widetilde\gamma\in\widetilde G(L)$ of $\gamma$.
Since $\gamma$ (and hence $\widetilde\gamma$) is semisimple, we have $Z_{\widetilde G}(\widetilde\gamma)=\pi^{-1}(Z_G(\gamma)^\circ)$.
Let $\widetilde\mu$ be the lift of $\mu$ with $\kappa(\vp^{\widetilde\mu})=\kappa(\widetilde\gamma)$, so that $Y_{\widetilde\mu}(\widetilde\gamma)\ne\varnothing$.
By \cite[Proposition 3.1]{HV18}, $\pi$ induces a universal homeomorphism from each connected component of the affine Grassmannian of $\widetilde G$ onto its image in $G(L)/K$.
Thus $\pi$ induces bijections $\widetilde T(L)\backslash\Sigma^{\tp}(Y_{\widetilde\mu}(\widetilde\gamma))\xrightarrow{\sim}\Sigma^{\tp}(Y_\mu(\gamma))$ and
$$
Z_{\widetilde G}(\widetilde\gamma)(L)\backslash\Sigma^{\tp}(Y_{\widetilde\mu}(\widetilde\gamma))
\xrightarrow{\sim} Z_G(\gamma)^\circ(L)\backslash\Sigma^{\tp}(Y_\mu(\gamma)).
$$
Thus $\sharp\bigl(Z_G(\gamma)^\circ(L)\backslash\Sigma^{\tp}(Y_\mu(\gamma))\bigr)=\dim V_{\widetilde\mu}^{\widehat{\widetilde G}}(\lambda_{\widetilde\gamma})=\dim V_\mu^{\widehat G}(\lambda_\gamma)$, as desired.
\end{proof}
\bibliographystyle{myamsplain}
\bibliography{reference}
\end{document}